\documentclass[12pt]{article}
\usepackage{amsmath,amsfonts,amsthm,amssymb,eucal}
\newcommand{\Qp}{\mathbb Q_p}
\newcommand{\Cp}{\mathbb C_p}
\newcommand{\Zp}{\mathbb Z_p}

\newcommand{\LA}{\mathcal L_A}
\newcommand{\LV}{\mathcal L_A}
\newcommand{\B}{\mathcal B}
\newcommand{\MA}{\mathcal M(\LA)}

\newcommand{\AAA}{\mathfrak A}
\newcommand{\UU}{\mathfrak U_1(K)}
\newcommand{\UUQ}{\mathfrak U_1(\Qp )}
\DeclareMathOperator{\Gal}{Gal}

\begin{document}
\newtheorem*{teo}{Theorem}

\pagestyle{plain}
\title{Non-Archimedean Unitary Operators}
\author{Anatoly N. Kochubei\\
\footnotesize Institute of Mathematics,\\
\footnotesize National Academy of Sciences of Ukraine,\\
\footnotesize Tereshchenkivska 3, Kiev, 01601 Ukraine
\\ \footnotesize E-mail: \ kochubei@i.com.ua}
\date{}
\maketitle

\bigskip
\begin{abstract}
We describe a subclass of the class of normal operators on Banach
spaces over non-Archimedean fields (A. N. Kochubei, J. Math. Phys. 51
(2010), article 023526) consisting of operators whose properties
resemble those of unitary operators. In particular, an analog of
Stone's theorem about one-parameter groups of unitary operators is
proved.
\end{abstract}


\bigskip
\section{INTRODUCTION}

{\bf 1.1}. In a previous paper \cite{K1}, we found a class of
non-Archimedean normal operators, bounded linear operators on Banach
spaces over non-Archimedean fields possessing orthogonal, in the
non-Archimedean sense, spectral decompositions. It is a natural problem
now to find out what operators in the non-Archimedean setting should be
seen as unitary ones. Classically, the correspondence between
selfadjoint and unitary operators extends, via the well-known
functional calculus, the correspondence $\lambda \mapsto e^{i\lambda }$
between real numbers and complex numbers from the unit circle.

There is no direct analog of the function $e^{i\lambda }$ in the
non-Archimedean case. However we will see that its natural counterpart
in the context of non-Archimedean operator theory is the function
$\lambda \mapsto z^{\lambda }$ where $\lambda$ runs the ring $\Zp$ of
$p$-adic integers, $z$ belongs to the group of principal units of a
non-Archimedean field $K$ (in another language, $z$ is a positive
element of $K$ \cite{Sch}). The image of this function also belongs to
the group of principal units. This prompts to define a non-Archimedean
unitary operator as an operator of the form $I+V$ where $I$ is the unit
operator, $V$ is a normal operator in the sense of \cite{K1}, $\|V\| <
1$. The normality assumption is essential -- otherwise $I+V$ can be
non-diagonalizable together with $V$. This shows (see also \cite{Shil})
that the isometricity is not a substitute of unitarity in the
non-Archimedean case. In fact we will use a refined version of the
above definition; see Section 2.

In classical operator theory, the main result about unitary operators
is Stone's theorem about the representation of a one-parameter unitary
group in the form $t\mapsto e^{itA}$ where $t\in \mathbb R$, $A$ is a
selfadjoint operator. We find its non-Archimedean analog -- a
one-parameter group parametrized by the group of principal units of
$\Qp$ has the form $s^A$ where $A$ belongs to the class of normal
operators in the sense of \cite{K1}, $\|A\| \le 1$. This result can be
reformulated from the setting with the parameter from the group of
principal units to the case of the parameter from $\Zp$.

{\bf 1.2.} Let us recall principal notions and results from \cite{K1}.
We will not explain the basic notions of non-Archimedean analysis; see
\cite{PS,R,vR,Sch}.

Let $A$ be a bounded linear operator on a Banach space $\B$ over a
complete non-Archimedean valued field $K$ with a nontrivial valuation;
$|\cdot |$ will denote the absolute value in $K$. Denote by $\LA$ the
commutative Banach algebra generated by the operators $A$ and $I$.
$\LA$ is a closure of the algebra $K[A]$ of polynomials in $A$, with
respect to the norm of operators; thus $\LA$ is a Banach subalgebra of
the algebra $L(\B )$ of all bounded linear operators. Elements $\lambda
\in K$ are identified with the operators $\lambda I$.

The spectrum $\MA$ is defined (see \cite{Berk}) as the set of all
bounded multiplicative seminorms on $\LA$. In a natural topology, it is
a nonempty Hausdorff compact topological space. If the algebra $\LA$ is
uniform, that is $\|T^2\|=\|T\|^2$ for any $T\in \LA$, and all its
characters take their values in $K$, then \cite{Berk} the space $\MA$
is totally disconnected and $\LA$ is isomorphic to the algebra $C(\MA
,K)$ of continuous functions on $\MA$ with values from $K$. In this
case, under the above isomorphism, the characteristic functions
$\eta_\Lambda$ of nonempty open-closed subsets $\Lambda \subset \MA$
correspond to idempotent operators $E(\Lambda )\in \LA$, $\|E(\Lambda
)\|=1$. These operators form a finitely additive norm-bounded
projection-valued measure on the algebra of open-closed sets, with the
non-Archimedean orthogonality property
$$
\|f\| =\sup\limits_\Lambda \|E(\Lambda )f\|,\quad f \in \B.
$$

An operator $A$, for which the above picture takes place, is called
{\it normal}. We will call a normal operator {\it strongly normal}, if
its spectrum $\sigma (A)$ is a nonempty totally disconnected compact
subset of $K$, and $\MA =\sigma (A)$. A strongly normal operator admits
the spectral decomposition
$$
A=\int\limits_{\sigma (A)}\lambda \,E(d\lambda ).
$$

More generally, we get a functional calculus assigning to any
$K$-valued continuous function $\varphi$ the operator
$$
\varphi (A)=\int\limits_{\sigma (A)}\varphi (\lambda )\,E(d\lambda ),
$$
such that
$$
\| \varphi (A)\| \le \sup\limits_{\lambda \in \sigma (A)}|\varphi
(\lambda )|.
$$

Some sufficient conditions of strong normality were found in \cite{K1}.
Let $\dim \B <\infty$, and $\B =K^n$, with the norm $\|(x_1,\ldots
,x_n)\|=\max\limits_{1\le i\le n}|x_i|$. An operator $A$ is
represented, with respect to the standard basis in $K^n$, by a matrix
$(a_{ij})_{i,j=1}^n$. Its operator norm coincides with $\|A\|
=\max\limits_{i,j}|a_{ij}|$ (see \cite{Se}). It is sufficient to
consider the case where $\|A\|=1$.

Let $\widehat{K}$ be the residue field of the field $K$. Together with
the operator $A$, we consider its reduction, the operator $\AAA$ on the
$\widehat{K}$-vector space $\widehat{\B }=\widehat{K}^n$ corresponding
to the matrix $(\widehat{a_{ij}})$ where $\widehat{a_{ij}}$ is the
image of $a_{ij}$ under the canonical mapping $O\to \widehat{K}$ ($O$
is the ring of integers of the field $K$). An operator $A$ is called
nondegenerate, if $\AAA \ne \nu I$ for any $\nu \in \widehat{K}$.

It was proved in \cite{K1} that $A$ is strongly normal, if it is
nondegenerate, all its eigenvalues belong to $K$, and its reduction
$\AAA$ is diagonalizable. These conditions are satisfied, for example,
if $\AAA$ has $n$ different eigenvalues from $\widehat{K}$.

In the infinite-dimensional situation, a similar result holds \cite{K1}
(with the representation of operators by infinite matrices), if we
assume in addition that $K$ is algebraically closed, $\B$ is the space
of sequences tending to zero, $A$ is a bounded operator with a compact
spectrum, and the resolvent of $A$ belongs, in a weak sense, to the
space of Krasner analytic functions outside the spectrum. For example,
if a compact operator (that is a norm limit of a sequence of finite
rank operators) is such that its reduction is diagonalizable, then it
is strongly normal.

Note that for a strongly normal operator $A$ and any continuous
$K$-valued function $\varphi$ on $\sigma (A)$, the operator $B=\varphi
(A)$ is strongly normal. Indeed, considering the functional model of
the algebra $\LA$ we see that the spectrum of the operator $B$
coincides with the set $f(\sigma (A))$. The Banach algebra $\mathcal
L_B$ is a subalgebra of $\LA$, and its functional model coincides with
the closure in $C(\sigma (A),K)$ of the set of functions $\pi \circ f$
where $\pi$ is an arbitrary polynomial. The convergence of the sequence
$\pi_n\circ f$ in $C(\sigma (A),K)$ is equivalent to the convergence of
the sequence of polynomials $\pi_n$ in the space $C(\sigma (B),K)$. By
Kaplansky's theorem (see Theorem 43.3 in \cite{Sch}), $\mathcal L_B$ is
isometrically isomorphic to $C(\sigma (B),K)$.

\section{Unitary Operators}

{\bf 2.1.} An operator $U$ on a Banach space $\B$ over a complete
non-Archimedean valued field $K$ with a nontrivial valuation will be
called {\it unitary}, if $U=I+V$ where $\|V\|<1$ and $V$ is strongly
normal. A unitary operator admits a spectral decomposition
$$
U=\int\limits_{\sigma (V)}(1+\lambda )E_V(d\lambda
)=\int\limits_{\sigma (U)}\mu E_U(d\mu ).
$$
Here $E_V$ is the spectral measure of the operator $V$, the mapping
$\varphi (\lambda )=1+\lambda$ transforms the spectrum of $V$ into that
of $U$, $E_U(M)=E_V(\varphi^{-1}(M))$ for any open-closed subset of
$\sigma (U)$.

Below we assume that the field $K$ is an extension of the field $\Qp$
of $p$-adic numbers, and the absolute value $|\cdot |$ is an extension
of the $p$-adic absolute value.

Denote by $\UU$ the group of principal units of the field $K$, that is
$\UU=\{ 1+\lambda :\ \lambda \in K,|\lambda |<1\}$. We will consider
{\it one-parameter groups} $U(s)$, $s\in \UU$, of unitary operators,
that is families of unitary operators, continuous with respect to the
norm of operators, such that
\begin{equation}
U(s_1s_2)=U(s_1)U(s_2),\quad s_1,s_2\in \UU .
\end{equation}

A one-parameter group of unitary operators can be constructed as
follows. Let $A$ be a strongly normal operator, $\|A\|\le 1$, $\sigma
(A)\subseteq \Zp$. Consider the function $f_s(\lambda )=(1+z)^\lambda$
where $s=1+z$, $z\in K$, $|z|<1$, $\lambda \in \Zp$. This function can
be defined by its Mahler expansion \cite{R,Sch}:
$$
f_s(\lambda )=\sum\limits_{n=0}^\infty z^nP_n(\lambda )
$$
where
$$
P_n(\lambda )=\frac{\lambda (\lambda -1)\cdots (\lambda -n+1)}{n!},\
n\ge 1;\quad P_0(\lambda )\equiv 1.
$$
An equivalent definition \cite{R,Sch} can be made in terms of the
approximation of a $p$-adic integer $\lambda$ by a sequence of
nonnegative integers, for which the function is defined in the
straightforward way.

Set
\begin{equation}
U(s)=(1+z)^A=\int\limits_{\sigma (A)}(1+z)^\lambda E_A(d\lambda ),\quad
s=1+z\in \UU.
\end{equation}
Due to the non-Archimedean orthonormality of the Mahler basis, we have
$$
(1+z)^\lambda =1+v_z(\lambda ),\quad v_z(\lambda
)=\sum\limits_{n=1}^\infty z^nP_n(\lambda ),
$$
$$
\sup\limits_{\lambda \in \Zp}\left| \sum\limits_{n=1}^\infty
z^nP_n(\lambda )\right| =|s|<1,
$$
so that $U(s)$ is a unitary operator, that is $U(s)=I+V(s)$ where
$$
V(s)=\int\limits_{\sigma (A)}v_z(\lambda )E_A(d\lambda
)=\int\limits_{\sigma (V(s))}\mu E_{V(s)}(d\mu ),\quad
E_{V(s)}(M)=E_A(v_z^{-1}(M)).
$$

It follows from the approximative description of the function $f_s$
that $U(s)$ possesses the required group property. Next, let
$s_1=1+z_1$, $s_2=1+z_2$, $|z_1|<1$, $|z_2|<1$. Using (2) we find that
$$
U(s_1)-U(s_2)=\int\limits_{\sigma (A)}\left[ \sum\limits_{n=1}^\infty
\left( z_1^n-z_2^n\right) P_n(\lambda )\right] E_A(d\lambda ),
$$
so that
$$
\|U(s_1)-U(s_2)\| \le \sup\limits_{\lambda \in \Zp}\left|
\sum\limits_{n=1}^\infty \left( z_1^n-z_2^n\right) P_n(\lambda )\right|
\le \sup\limits_{n\ge 1}\left| z_1^n-z_2^n\right| =|s_1-s_2|.
$$
Therefore the function $s\mapsto U(s)$ is continuous with respect to
the operator norm.

\medskip
{\bf 2.2.} In particular, the above construction makes sense for $s\in
\UUQ$, and the formula (2) defines a norm-continuous one-parameter
group of unitary operators $\UUQ\ni s\mapsto U(s)$. The next result is
a converse statement, an analog of Stone's theorem.

\medskip
\begin{teo}
Let $U(s)$, $s\in \UUQ$, $p\ne 2$, be a norm continuous one-parameter
group of unitary operators, such that the spectrum of the strongly
normal operator $U(1+p)-I$ is contained in $p\Zp$. Then there exist
such a strongly normal operator $A$, $\sigma (A)\subseteq \Zp$, that
$U(s)=s^A$, $s\in \UUQ$.
\end{teo}

\medskip
{\it Proof}. Each element $s\in \UUQ$ can be represented, in a unique
way, as
\begin{equation}
s=(1+p)^\zeta ,\quad \zeta \in \Zp .
\end{equation}

Indeed, set $\zeta =\dfrac{\log s}{\log (1+p)}$ (see \cite{R,Sch}
regarding properties of the $p$-adic logarithm). We have $\zeta \in
\Qp$, $|\log s|\le p^{-1}$, $|\log (1+p)|=|p|=p^{-1}$, so that $\zeta
\in \Zp$. On the other hand, $\exp (\zeta \log (1+p))=(1+p)^\zeta$
(\cite{Sch}, Theorem 47.10), which implies (3).

Let us write the canonical representation
$$
\zeta =\zeta_0+\zeta_1p+\zeta_2p^2+\cdots ,\quad \zeta_j\in
\{0,1,\ldots ,p-1\}.
$$
The series $(1+p)^\zeta =\sum\limits_{n=0}^\infty p^nP_n(\zeta )$
converges uniformly with respect to $\zeta \in \Zp$, so that the
function $\zeta \mapsto (1+p)^\zeta$ is continuous. Due to the norm
continuity of $U$,
\begin{equation}
U(s)=\lim\limits_{n\to \infty}[U(1+p)]^{\zeta_0+\zeta_1p+\cdots
+\zeta_np^n}.
\end{equation}

Denote
$$
a_n(\lambda )=(1+p)^{(\zeta_0+\zeta_1p+\cdots +\zeta_np^n)\lambda
},\quad \lambda \in \Zp .
$$
Let us prove that
\begin{equation}
a_n(\lambda )\longrightarrow (1+p)^{\zeta \lambda },\quad \text{as
$n\to \infty$},
\end{equation}
uniformly with respect to $\lambda$.

Indeed, we have the estimate
\begin{multline*}
\left| a_n(\lambda )-(1+p)^{\zeta \lambda }\right| =\left|
(1+p)^{(\zeta_{n+1}p^{n+1}+\zeta_{n+2}p^{n+2}+\cdots )\lambda
}-1\right|
\\
=\left| \sum\limits_{k=1}^\infty
p^kP_k((\zeta_{n+1}p^{n+1}+\zeta_{n+2}p^{n+2}+\cdots )\lambda )\right|
\le \sup\limits_{k\ge 1}p^{-k}\left|
P_k((\zeta_{n+1}p^{n+1}+\zeta_{n+2}p^{n+2}+\cdots )\lambda )\right|
\end{multline*}
where
$$
\left| P_k((\zeta_{n+1}p^{n+1}+\zeta_{n+2}p^{n+2}+\cdots )\lambda
)\right| \le p^{-n-1}|k!|^{-1}\le p^{-n-1+\frac{k-1}{p-1}},
$$
so that, uniformly with respect to $\lambda \in \Zp$,
$$
\left| a_n(\lambda )-(1+p)^{\zeta \lambda }\right| \le
p^{-n-1}\sup\limits_{k\ge 1}p^{-k+\frac{k-1}{p-1}}\longrightarrow 0,
$$
as $n\to \infty$.

Now we return to the expression (4). By our assumption, $U(1+p)=I+V$
where $V$ is a strongly normal operator, $\sigma (V)\subseteq p\Zp$. We
have
$$
U(1+p)=\int\limits_{\sigma (V)}(1+\lambda )E_V(d\lambda )
$$

The Banach algebra $\LV$ generated by the operator $V$ contains the
strongly normal operator
$$
A=\frac1{\log (1+p)}\log (I+V)=\frac1{\log
(1+p)}\sum\limits_{k=1}^\infty \frac{(-1)^{k-1}}{k}V^k
=\int\limits_{\sigma (V)}\frac{\log (1+\lambda )}{\log
(1+p)}E_V(d\lambda ).
$$
Obviously, $\sigma (A)\subseteq \Zp$ and
$$
(1+p)^{\frac{\log (1+\lambda )}{\log (1+p)}}=1+\lambda ,
$$
so that $U(1+p)=(1+p)^A$, and it follows from (4) that
$$
U(s)=\lim\limits_{n\to \infty}\left[
(1+p)^A\right]^{\zeta_0+\zeta_1p+\cdots +\zeta_np^n}.
$$
Switching to the functional model and using (5) and (3) we obtain the
required formula for the operators $U(s)$. $\qquad \blacksquare$

\medskip
Note that the condition regarding the operator $U(1+p)-I$ is satisfied
automatically, if $U(1+p)=I+V$, $\|V\|<1$, and $K=\Qp$.

\medskip
{\bf 2.3.} Let $W(z)$, $z\in \Zp$, $p\ne 2$, be a norm continuous
unitary representation of the additive group $\Zp$, and the spectrum of
the operator $W(p^{-1}\log (1+p))-I$ lies in $p\Zp$. Denote $s=e^{pz}$,
$U(s)=W(z)$. Then $s\in \UUQ$, and $s\mapsto U(s)$ is a one-parameter
group satisfying the conditions of the above Theorem. We obtain the
expression
$$
W(z)=e^{pzA},\quad z\in \Zp ,
$$
where $A$ is a strongly normal operator, $\sigma (A)\subseteq \Zp$.

\section{Example. Galois Representations}

In this section we follow \cite{Berg,Fon,Sen}.

Let $K$ be a finite extension of $\Qp$, and $\varepsilon^{(n)}\in
\bar{K}$ ($\bar{K}$ is an algebraic closure of $K$) is a sequence of
primitive roots of unity of orders $p^n$, such that
$$
\varepsilon^{(0)}=1,\quad \varepsilon^{(1)}\ne 1,\quad
(\varepsilon^{(n+1)})^p=\varepsilon^{(n)},\ n=0,1,\ldots.
$$

Denote $K_n=K(\varepsilon^{(n)})$, $K_\infty
=\bigcup\limits_{n=0}^\infty K_n$, $G_K=\Gal (\bar{K}/K)$. Let
$\mu_{p^n}$ be the set of roots of unity of order $p^n$; thus
$\varepsilon^{(n)}\in \mu_{p^n}$, $n\ge 0$. The cyclotomic character
$\chi:\ G_K\to \Zp^*=\{ x\in \Zp :\ |x|=1\}$ is defined via the
equality
$$
\sigma (\zeta )=\zeta^{\chi (\sigma )},\quad \text{for all $\sigma \in
G_K$, $\zeta \in \mu_{p^\infty}=\bigcup\limits_{n=0}^\infty
\mu_{p^n}$.}
$$
$\chi$ is continuous with respect to the standard topology of $G_K$ as
a profinite group.

The kernel of the cyclotomic character coincides with $H_K=\Gal
(\bar{K}/K_\infty )$. Therefore $\chi$ identifies $\Gamma_K=\Gal
(K_\infty /K)=G_K/H_K$ with an open subgroup of the multiplicative
group $\Zp^*$.

By definition, {\it a $p$-adic representation} $V$ of the group $G_K$
is a finite-dimensional vector space over $\Qp$ with a continuous
linear action of $G_K$.

Let $\widehat{K}_\infty$ be the $p$-adic completion of $K_\infty$. Let
us consider the action of $\Gamma_K$ on the $\widehat{K}_\infty$-vector
space $\left( \Cp \otimes_{\Qp}V\right)^{H_K}$ of elements from $\Cp
\otimes_{\Qp}V$ fixed under the action of $H_K$. If $d=\dim_{\Qp}V$,
then $\left( \Cp \otimes_{\Qp}V\right)^{H_K}$ is a
$\widehat{K}_\infty$-vector space of dimension $d$. The group
$\Gamma_K$ acts on the union $\mathbb D_{\text{Sen}}(V)$ of
finite-dimensional subspaces of $\left( \Cp
\otimes_{\Qp}V\right)^{H_K}$ invariant with respect to $\Gamma_K$, and
$\dim_{K_\infty}\mathbb D_{\text{Sen}}(V)=d$.

By Sen's theorem \cite{Sen}, there is a unique $K_\infty$-linear
operator $\Theta_V$ on $\mathbb D_{\text{Sen}}(V)$, such that for any
$\omega \in \mathbb D_{\text{Sen}}(V)$, there exists such an open
subgroup $\Gamma_\omega \subset \Gamma_K$ that
\begin{equation}
\sigma (\omega )=\left[ \exp \left( \Theta_V\log \chi (\sigma )\right)
\right] \omega
\end{equation}
for all $\sigma \in \Gamma_\omega$.

A representation $V$ is called a {\it Hodge-Tate representation} if,
for a certain basis $e_1,\ldots ,e_d\in \mathbb D_{\text{Sen}}(V)$, the
operator $\Theta_V$ is diagonal, with eigenvalues from $\mathbb Z$. In
this case, we can introduce a norm in $\mathbb D_{\text{Sen}}(V)$
setting
$$
\|x_1e_1+\cdots +x_de_d\|=\max (|x_1|,\ldots ,|x_d|),\quad x_1,\ldots
,x_d\in K_\infty .
$$
Then $\Theta_V$ is obviously strongly normal, $\|\Theta_V\| \le 1$.

Taking into account the continuity of $\chi$ we can choose a so small
open subgroup $\Lambda \subset \Gamma_K$ that $|\chi (\sigma )-1| \le
p^{-1}$ for all $\sigma \in \Lambda$. For every $\sigma \in \Lambda$,
the right-hand side of (6) defines a unitary operator on $\mathbb
D_{\text{Sen}}(V)$.

\section*{ACKNOWLEDGEMENT}
This work was supported in part by the Ukrainian Foundation
for Fundamental Research, Grant 29.1/003.

\end{document}